\documentclass[12pt]{amsart}
\usepackage{amssymb}
\usepackage[T1]{fontenc}
\usepackage{amsmath,stackrel}
\usepackage{indentfirst}
\usepackage{xspace}
\usepackage{color}
\usepackage{fancyhdr} 
\usepackage{fancybox}
\usepackage[brazil,english]{babel}
\usepackage[latin1]{inputenc}
\usepackage{graphicx}
\usepackage[all]{xy}
\usepackage{hyperref}

\theoremstyle{plain}
\newtheorem{teo}{Theorem}

\newtheorem{prop}[teo]{Proposition}

\theoremstyle{definition}
\newtheorem{defi}{Definition}[section]
\newtheorem{exe}{Example}

\theoremstyle{remark}
\newtheorem*{obs}{Remark}
\numberwithin{teo}{section}

\newcommand{\C}{\mathbb{C}} 

\newenvironment{prova}
{{\em{\noindent \textbf{Proof:} }}
}
{\hfill $\blacksquare$}

\DeclareMathOperator{\rank}{rank}

\DeclareMathOperator{\codim}{codim}

\begin{document}

\title{Generic sections of  essentially isolated determinantal singularities}
\author{Jean-Paul Brasselet, Nancy Chachapoyas and Maria A. S. Ruas}
\date{}
\thanks{ The first author is partially supported by CNPq grant no. 400580/2012-8. The second and third
		authors are partially supported by FAPESP grants no. 2010/09736-1, 2011/20082-6 and 2014/00304-2.
		The third author is partially supported by CNPq grant no. 301474/2005-2.}
\begin{abstract}
  We study the essentially isolated determinantal singularities (EIDS), defined  by W. Ebeling and S. Gusein-Zade \cite{Gusein}, as a generalization of isolated singularity. We prove in dimension $3$ a minimality theorem for the Milnor number of a generic hyperplane section of an EIDS, generalizing previous results by J. Snoussi in dimension $2$. We define strongly generic hyperplane sections of an EIDS and show that  they are still EIDS.
 Using strongly general hyperplanes, we extend a result of L\^e D. T. concerning constancy of the Milnor number.
\end{abstract}

\maketitle

\thispagestyle{plain}
\pagestyle{fancy}
\fancyhead{}
\fancyhead[CO]{\bfseries \small{Generic sections of  essentially isolated determinantal singularities}}
\fancyhead[LE]{\small{Jean-Paul Brasselet, Nancy Chachapoyas and Maria A. S. Ruas } }
\renewcommand{\headrulewidth}{0pt}


\section{Introduction}
\label{sec:introduction}

\thispagestyle{empty}
In this work, we study the essentially isolated determinantal singularities (EIDS), which have been defined by W.~{E}beling and S.~M. Guse{\u\i}n-Zade in \cite{Gusein}. This type of singularities is a natural generalization of isolated ones. A generic determinantal variety $M_{m,n}^t$ is the subset of the space of $m \times n$ matrices, given by the matrices of rank less than 
$t$, where $t\leqslant \min\{m,n\}$.
 A variety $X\subset \C^N$ is  determinantal if $X$ is the pre-image of $M_{m,n}^t$ by a holomorphic function $F:\C^N\to M_{m,n}$ with the condition that the codimension of $X$ in $\C^N$ is the same as the codimension of $M_{m,n}^t$ in $M_{m,n}$.

Determinantal varieties have isolated singularities if $N\leqslant (m-t+2)(n-t+2)$  and they admit smoothing if $N<(m-t+2)(n-t+2)$. Several recent works investigate determinantal varieties with isolated singularities. The Milnor number of a determinantal surface was defined in \cite{Pike, Bruna, MIRIAMSP} while the vanishing Euler characteristic of a determinantal variety was defined in \cite{Pike, Bruna}. Other recent results on isolated determinantal varieties related to this paper  appear in particular in \cite{Anne-Matthias, Gaffney-Antoni}.

In this work we study  the set of limits of tangent hyperplanes to determinantal varieties. For determinantal surfaces in $\C^4$ and $3$-fold singularities in $\C^5$,  we charac\-terize these hyperplanes by the fact that the Milnor number of their sections with the surface in the first case or the $3$-dimensional  determinantal variety in the second case, is not minimum. The first case follows as a consequence of results of J. Snoussi  in \cite{Jawad}.

We also prove that  given a $d$- dimensional EIDS $X$, given $H$ and $H'$ strongly general hyperplanes to $X$ (Definition \ref{strongly}), there are $d-2$ linear planes $P\subset H$ and $P'\subset H'$ contained in $H$ and $H'$, such that the Milnor number of the surfaces $X\cap P$ and $X\cap P'$ are equal. In the case that the generic section is a curve the result has been already proved by L\^e D. T. in \cite{Le}.




\section{Essentially isolated determinantal singularity }

In this section we give the definition and basic results on essentially isolated determinant singularities, following \cite{Gusein}.

We denote by $M_{m,n}$ the set of matrices $m\times n$ with complex entries.
\begin{defi} For all  $t$, $ 1\leqslant t\leqslant \min\{m,n\}$, let $M_{m,n}^{t}$ be  the subset of $M_{m,n}$ whose elements are matrices of rank less than $t$:
$$M_{m,n}^{t}=\{A\in M_{m,n} \vert  \rank(A)< t \}.$$
This set is a singular variety of codimension $ (m-t+1)(n-t+1) $ in $M_{m,n}$, called generic determinantal variety.
\end{defi}
The singular set of $M_{m,n}^{t}$ is $M_{m,n}^{t-1}$.
The partition of $M_{m,n}^{t}$ defined by 
$$M_{m,n}^{t}=\cup_{i=1,...,t}( M_{m,n}^{i} \backslash M_{m,n}^{i-1})$$ is a Whitney stratification \cite{Harris}.

Let $F: \C^{N}\rightarrow M_{m,n}$ be a map defined by $F(x)=(f_{ij}(x))$, whose entries are complex analytic functions defined on an open domain  $U\subset \C^N$.

\begin{defi}

 The analytic variety $X=F^{-1}(M_{m,n}^{t})$ is called determinantal variety of type $(m,n,t)$, if $\codim X=\codim M_{m,n}^t=(m-t+1)(n-t+1)$.
\end{defi}

A generic map $F$ intersects transversally the strata $M_{m,n}^i \backslash M_{m,n}^{i-1}$ of the variety $M_{m,n}^t$.
\begin{defi} A point $x\in X= F^{-1}(M_{m,n}^{t})$ is called essentially nonsingular if at this point the map $F$ is transversal to the corres\-ponding stratum of the variety $M_{m,n}^{t}$ (that is, to $\{M_{m,n}^{i} \backslash M_{m,n}^{i-1}\}$, where $i=\rank F(x)+1$ ).
\end{defi}

%
%
%

\begin{defi}\label{EIDS} A germ $(X,0)\subset (\C^{N},0) $ of a determinantal variety of type $(m,n,t)$ has an essentially isolated singular point at the origin (or is an essentially isolated determinantal singularity, EIDS) if it has only essentially nonsingular points in a punctured neighborhood of the origin in $X$.
\end{defi}

An EIDS $X\subset \C^N$ has isolated singularity  if and only if $N\leqslant (m-t+2)(n-t+2)$. An EIDS with isolated singularity will be called isolated determinantal singularity, denoted by IDS.

We want to consider deformations of an EIDS that are themselves determinantal varieties of the same type.

\begin{defi} An essential smoothing $\widetilde{X}_s$ of the EIDS $(X,0)$ is a subvariety lying in a neighborhood $U$ of the origin in $\C^N$ and defined by $\widetilde{X}_s = \widetilde{F}^{-1}_s (M_{m,n}^t)$ where 
 $\widetilde{F}:U\times \C \to M_{m,n}$ is a perturbation of the germ $F$, with $\widetilde{F}_{s}(x)=\widetilde{F}(x,s)$, $\widetilde{F}_0(x)=F(x)$ and such that $\widetilde{F}_s:U\to M_{m,n}$ is transversal to all strata $M_{m,n}^i \backslash M_{m,n}^{i-1}$.
\end{defi}
An essential smoothing is in general not smooth (when $ N\geq (m-t+2)(n-t+2)$) as we see in the following theorem.

\begin{teo}\cite{Wahl} \label{smoothing} Let $(X,0)\subset (\C^N,0)$ be the germ of a determinantal variety with isolated singularity at the origin. Then, $X$ has a smoothing if and only if  $N<(m-t+2)(n-t+2).$
\end{teo}

The singular set of the essential smoothing $\widetilde{X}_s$ is $\widetilde{F}_s^{-1}(M^{t-1}_{m,n})$. Since $\widetilde{F}$ is transversal to the strata of the Whitney stratification $M_{m,n}^t$, the partition  $\widetilde{X}_s=\cup_{1\leq i \leq t} \widetilde{F}_s^{-1}(M^{i}_{m,n}\backslash M^{i-1}_{m,n})$ is a Whit\-ney stratification of $\widetilde{X}_s$.


\begin{exe} Let $X=F^{-1}(M_{2,3}^2)$, where
 $$\begin{array}{cccl}
                F :& \C^{4}      & \rightarrow & M_{2,3} \\
                   & (x,y,z,w)   & \mapsto     &\left(
                                                 \begin{array}{ccc}
                                                   z & y+w & x \\
                                                   w & x   & y \\
                                                 \end{array}
                                               \right)

              \end{array}
$$
The following matrix defines an essential smoothing $\widetilde{X}_s=\widetilde{F}_s^{-1}(M_{2,3}^2)$ of $X$
$$\begin{array}{cccl}
                \widetilde{F} :& \C^{4}  \times \C   & \rightarrow & M_{2,3} \\
                   & (x,y,z,w,s)   & \mapsto     &\left(
    \begin{array}{ccc}
      z & y+w & x+s \\
      w & x & y \\
    \end{array}
  \right)
              \end{array}$$
In this case the essential smoothing $\widetilde{X}_s$ is a smoothing.

\end{exe}

\section{L\^e-Greuel formula type for IDS with smoo\-thing}

We review in this section some results about isolated singularities admitting smoothing following \cite{MIRIAMSP}.

Let $(X,0)\subset (\C^N,0)$ be the germ of a $d$-dimensional variety with isolated singularity at the origin. Suppose that $X$ has a smoothing. Then, there exists a flat family $\pi: \widetilde{X}\subset U\times \C\to \C$ such that the fiber $X_s=\pi^{-1}(s)$ is smooth for all $s\neq 0$ and $X_0=X$.

Let $p:(X,0)\to \C$ be a complex analytic function defined in $X$ with isolated singularity at the origin. Let us consider a function
$$\begin{array}{cccc}
      \widetilde{p}:&   \C^N\times \C   &   \to &  \C \\
      & (x,s)  & \mapsto & \widetilde{p}(x,s),
\end{array}$$
such that $ \widetilde{p}(x,0)=p(x)$ and for all $s\neq 0$, $ \widetilde{p}(\cdot,s)=p_s$ is a Morse function on $X_s.$

Thus we have the following diagram

\begin{displaymath}
\xymatrix{ X_s \ar[d]^{p_s} \ar@{^{(}->}[r]  &  \C^N\times \C \ar[d]^{(\pi,\widetilde{p})}\\
                    \C\times \{s\}    \ar@{^{(}->}[r]     & \C\times \C}
\end{displaymath}

\begin{prop} \cite[Proposition 4.1]{MIRIAMSP}  \label{prop:Miriam} Let $X$ be a d-dimensional variety with isolated singularity at the origin admitting smoothing and $p_s:X_s \to \C$, $p_s=\widetilde{p}(\cdot,s)$ as above.
Then,

\begin{itemize}
\item [(a)] If $s\neq0$, $X_s\simeq p_s{-1}(0)\dot{\cup}\{ \text{cells of dimension d } \},$

\item [(b)] $\chi(X_s)=\chi(p_s^{-1}(0))+(-1)^d n_{0},$
\end{itemize}
where  $n_{0}$ is the number of critical points of $p_s$ and $\chi(X_s)$ denotes the Euler characteristic of $X_s$.
\end{prop}

The invariant $n_{0}$ is related to the polar multiplicity of $X$, $m_d(X)$ (\cite{MIRIAMSP}, see also \cite{Gaffney1}), in the following way:
\begin{defi}(The $d$-Polar multiplicity) Let ${X}$, $\widetilde{X}$, $p$ and $\widetilde{p}$ as above.
Let $P_d(X,\pi,{p})=\overline{\Sigma(\pi,\widetilde{p})|}_{\widetilde{X}_{reg}}$ be the relative polar variety of $X$ related to $\pi$ and $p$. We define $m_d(X,\pi,{p})=m_0(P_d(X,\pi,{p}))$.
\end{defi}

In general $m_d(X,\pi,{p})$ depends on the choices of $\widetilde{X}$ and $\widetilde{p}$. When the variety $X$ has a 
unique smoothing $\widetilde{X}$, then 
$m_d(X,\pi,{p})$ depends only on $X$ and ${p}$. If ${p}$ is a generic linear embedding, $m_d(X,p)$  is an invariant of the EIDS $X$, denoted by $m_d(X)$.

\begin{prop} \cite{MIRIAMSP} Under the conditions of Proposition \ref{prop:Miriam}, $n_{0}=m_d(X)$.
\end{prop}

%
%

\begin{teo}  \cite{greuellivro} Let $X_s$ be a smoothing of a normal isolated singularity, then $b_1(X_s)=0.$
\end{teo}

Let $X$ be a determinantal variety of type $(m,n,t)$ in $\C^N$ with $N< (m-t+2)(n-t+2)$. Then $X$ has isolated singularity and admits smoothing. 
\begin{defi} \cite{MIRIAMSP}\label{Defi:Milnor} Let $X$ be a determinantal surface in $\C^N$, with isolated singularity at the origin. The Milnor number of $X$, denoted by $\mu(X) $, is defined as the second Betti number of the generic fiber $X_s$,
$$\mu(X)=b_2(X_s).$$
\end{defi}


The following result appears in \cite{Pike,Bruna, MIRIAMSP}, for determinantal surfaces $X\subset \C^4$, but it also holds for any determinantal surface with isolated singularity in $\C^N$ admitting smoothing.

\begin{prop}\label{LeGreuel} Let $(X,0)\subset (\C^N,0)$ be the germ of a determinantal surface in $\C^N$ with isolated singularity at the origin admitting smoothing. Let $p:(\C^N,0)\to (\C,0)$  be a linear function whose restriction to $X$ has an isolated singularity at the origin. Then one has the L\^e-Greuel formula
\begin{equation}
\label{formulaGreuel}
\mu(X)+\mu(X\cap p^{-1}(0))=m_2(X).
\end{equation}
\end{prop}

When $d=\dim X> 2,$ the Betti numbers $b_i(X),$ $2 \leqslant i < d$ are not ne\-ces\-sarily zero (see \cite{Anne-Matthias}). In \cite{Pike, Bruna} the authors define the vanishing Euler charac\-teristic of varieties admitting smoothing.

\begin{defi} \cite{Bruna} Let $(X,0)\subset (\C^N,0)$ be an IDS such that $N < (m-t+2)(n-t+2)$. The vanishing Euler characteristic is defined by \[\nu(X)= (-1)^d(\chi(X_s)-1),\]
where $X_s$ is a smoothing of  $X$ and $\chi(X_s)$ is the Euler characteristic of $X_s$ .
\end{defi}

 \begin{teo} \label{Greuel} \cite{Bruna} Let $(X,0)\subset (\C^N,0)$ be an IDS such that $N<(m-t+2)(n-t+2)$ and let $p:(\C^N,0)\to (\C,0)$ be a linear projection whose restriction to $X$ has isolated singularity at the origin. Then,
\begin{equation}
\label{}
\nu(X)+\nu(X\cap p^{-1}(0))=m_d(X)
\end{equation}
\end{teo}

\begin{obs} When $d=2$, then $\nu(X)=\mu(X)$.
\end{obs}

\begin{exe} \cite{MP} Let $X=F^{-1}(M_{2,3}^2)\subset \C^4$ be the variety defined by:
$$\begin{array}{cccc}
      F:&  \C^4& \to &M_{2,3}\\
      &(x,y,z,w)   & \mapsto & \begin{pmatrix}
    z  & y+w  & x \\
    w  &  x & y
\end{pmatrix}
\end{array}.$$

The following matrix defines an essential smoothing $\widetilde{X}_s=\widetilde{F}_s^{-1}(M_{2,3}^2)$ of $X$
$$\begin{array}{cccl}
                \widetilde{F} :& \C^{4}  \times \C   & \rightarrow & M_{2,3} \\
                   & (x,y,z,w,s)   & \mapsto     &\left(
    \begin{array}{ccc}
      z & y+w & x+s \\
      w & x & y \\
    \end{array}
  \right)
              \end{array}$$
              
Consider $p:\C^4\to \C$ given by $p(x,y,z,w)=w$, then it follows that $m_2(X)=3$ and $\mu(X\cap p^{-1}(0))=2$, then $\mu(X)=1.$
\end{exe}

\section{General and strongly general hyperplanes}

The main goal of this section is to define general and strongly general hyperplanes over an EIDS in order to extend previous results by J. Snoussi and L\^e D. T. in the following sections. The two definitions are equivalent when the variety has an isolated singularity. The sections defined by the intersection of the IDS (respectively EIDS) and the general hyperplane (respectively strongly general hyperplane) determine another IDS (respectively EIDS).

\begin{defi} \label{Geral} The hyperplane $H\subset \C^N$, given by the kernel of the linear function $p:\C^N\to \C$ is called general to $X$ at $0$ if $H$ is transversal to all limits $T$ of tangent spaces to the regular part of $X$.
\end{defi}

\begin{prop} \cite{Jawad} A hyperplane $H$  is limit of tangent hyperplanes to $X$ at 0, if only if $H$ is not a general hyperplane.
\end{prop}

\begin{figure}
  \centering
  \includegraphics[width=12cm,height=7cm]{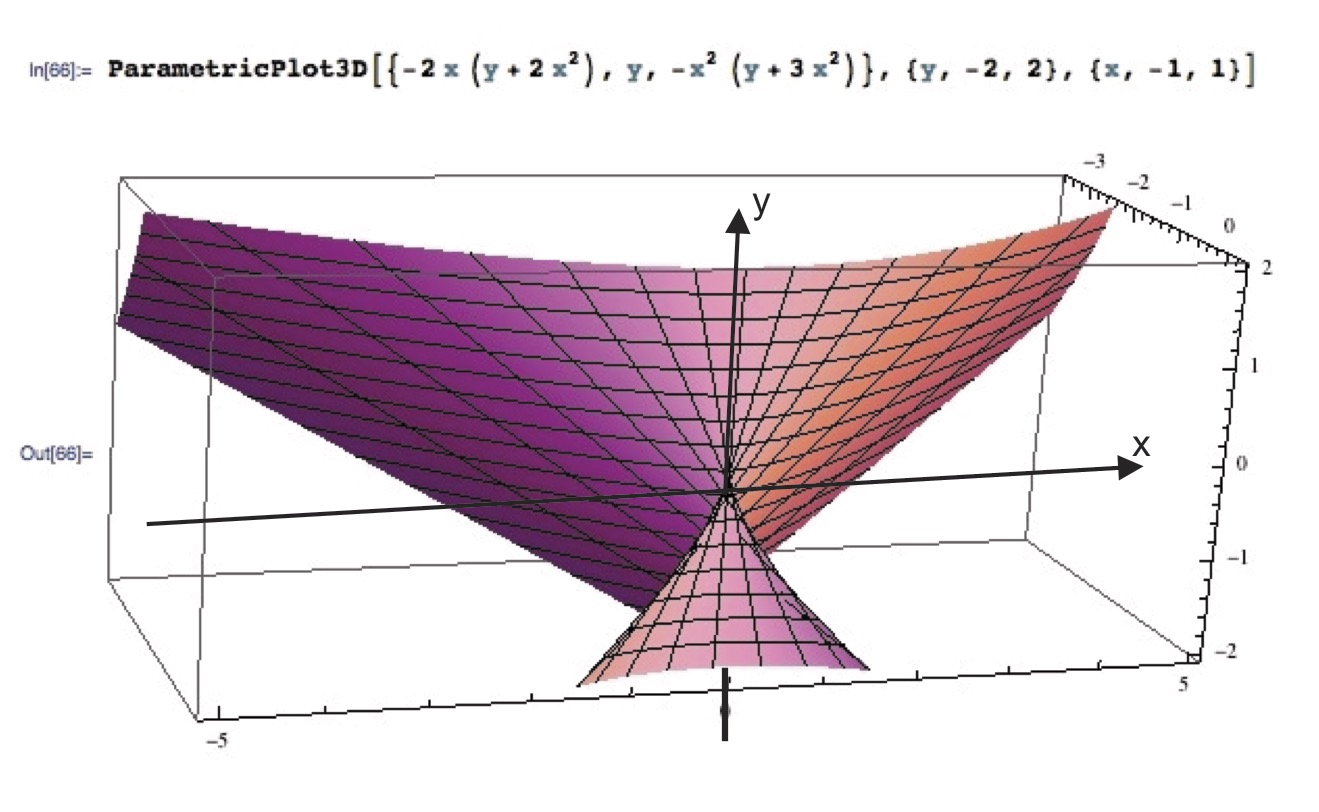}\\
  \caption{Swallowtail surface}\label{rabo}
\end{figure}

\begin{exe}
 Let $X$ be the swallowtail surface, given in  Figure $\ref{rabo}$. $X$ is the surface in $\C^3$ defined by the zeros of  $$256z^3-27x^4-128z^2y^2+144zx^2y+16zy^4-4x^2y^3=0.$$

The set of the limits of tangent hyperplanes to $X$ at zero is given by the hyperplane $z=0$ (see for example \cite{MMR}).

Any hyperplane different of the hyperplane $z=0$ is general, in particular the plane $H=\{(x,y,z): x=0\}$ is general, but $H$ is not transversal to the limits of lines tangent to the strata of dimension $1$.
\end{exe}

This example motivates the following definition:

\begin{defi}\label{strongly} Let $X\subset \C^N$ be a d- dimensional analytic  complex variety, and let  $\{V_{\lambda}\}_{\lambda\in\Lambda}$ be a stratification of $X$. The hyperplane $H\subset \C^N$ is called strongly general at the origin if it is general and there exists a neighborhood $U$ of $0$ such that for all strata $V_{\lambda}$ of $X$, with $0\in \overline{V}_{\lambda}$, we have that $H$ is transverse $V_{\lambda}$ at $x$, for all $x\in U\backslash \{0\}$.
\end{defi}

\begin{teo} \label{EIDS} Let $(X,0)\subset (\C^N,0)$ be an EIDS of type $(n,m,t)$. If $H\subset \C^N$ is a strongly general hyperplane at the origin, then $X\cap H$ is a $(d-1)$-dimensional EIDS in $\C^{N-1}$ of the same type.
\end{teo}
\begin{prova} Let $F:\C^N \to M_{m,n}$ be a function defining $X=F^{-1}(M_{m,n}^t)$. As $X$ is an EIDS, then $F$ is transversal to the strata of $M_{m,n}^t$ in $U\backslash \{0\}$, where $U$ is a sufficiently small neighborhood  of the origin. By hypothesis $H$ is transversal to the strata of $X$ outside the origin. Let $i:\C^{N-1}\to \C^N$ be a linear embedding such $i(\C^{N-1})=H$.
Then $F\circ i:\C^{N-1} \to M_{m,n}$ is transversal to all the strata of $M_{m,n}^t$ outside the origin, so $X\cap H=(F\circ i)^{-1}(M_{m,n}^t)\subset \C^{N-1}$ is an EIDS of the same type of $X$ and $\dim(X\cap H)=d-1$ in $\C^{N-1}.$
\end{prova}

\section{Minimality of the Milnor number}

The  minimality of the Milnor number of generic sections of  hypersurfaces with isolated singularities was studied by B. Teissier
\cite{Teissier3} and J.-P. Henry and L\^e D. T. \cite{Henry}.  T. Gaffney \cite{Gaffney} proved the result for  ICIS and J. Snoussi considered the case of  normal surfaces in $\C^N$.

In the following we denote by $|Z|$ the reduction of the variety $Z$ and   by $C_{X,x}$ the tangent cone to $X$ at $x$ (see \cite{Jawad}).

\begin{prop}\label{minimalidade sup} \cite[Theorem 4.2]{Jawad} Let $X\subset \C^N$ be the germ of a normal analytic surface. A hyperplane $H$ in $\C^N$ that does not contain any irreductible component of the tangent cone $|C_{X,0}| $   is general to $X$ at $0$ if and only if the section $X\cap H$ is reduced and the Milnor number $\mu(X\cap H)$ is minimum.
\end{prop}

\begin{prop} \cite{Le} \label{Legeral} Let $X\subset \C^N$ be a complex analytic variety of dimension $d$  and let $L$ be an affine subspace of $\C^N$ of codimension equal to the dimension of $C_{X,x}$ in $\C^N$ such that  $L\cap C_{X,x}=\{x\}$. Consider the projection $\pi:\C^N\to \C^{d}$ such that $\ker \pi=L$. Consider $\Delta _L$ the discriminant of $\pi$ restricted to $X$. Let $H_0$ be a hyperplane of $\C^{d}$ which is not a limit of tangent hyperplanes to the discriminant $\Delta _L$. Then, the hyperplane $H$ given by the inverse image $H=\pi^{-1}(H_0)$ is not  limit of tangent hyperplanes to $X$.
\end{prop}

\begin{prop}\label{Tessier} \cite[Corollary 2.3.2.1]{LeTessier} Let $X\subset \C^N$ be a reduced analytic variety, $c$ an integer and $E$ be a vector space in $\C^N$ of codimension $c$ given by the intersection of $c$ hyperplanes  $H_1, H_2, \dots, H_c$ such that each $H_i$ is not tangent to $|X\cap H_1\cap H_2\cap \dots \cap H_{i-1}|$, for $1 \leqslant i\leqslant c$. Then the multiplicity at $x$ of the polar variety $P_k(|X\cap E|)$ is equal to the multiplicity at $x$ of the polar varieties $P_k(X)$ for $0\leqslant k \leqslant d-c-1.$
\end{prop}

We extend Proposition \ref{minimalidade sup} to $3$-determinantal varieties with isolated singularities.

\begin{teo}Let $(X,0) \subset (\C^N,0)$ be the germ of a $3$- dimensional determinantal varie\-ty with isolated singularity 
and $H$ a hyperplane in $\C^N$. Suppose that  $X\cap H$ has an isolated singular point at the origin, then the following conditions are equivalent.
\begin{enumerate}
\item [({\rm i})] $H$ is general to $X$ at $0$.
\item [({\rm ii})] $\mu(X\cap H)$ is minimum and $\mu(X\cap H\cap H')$ is minimum for all $H'$ general to $X$ and  to $X\cap H$.
\end{enumerate}
\end{teo}

\begin{prova} Suppose that  $p:\C^N\to \C$ is general to  $X$ at the origin, $H=p^{-1}(0)$. Then $X\cap H$ is a determinantal surface with isolated singularity at the origin. Then there exists a Zariski open set $\Omega\subset \mathbb{P}^{N-1}$ such that $\forall$ $H'\in \Omega$, with $H'=\ker p'$,  $p':\C^N \to \C$ a linear map, we have that:
\begin{enumerate}
\item $ p$ is general to $X$ at the origin and $p$ is general to $X\cap H'$ at the origin.
\item $\mu(X\cap H')$ is minimum.
\end{enumerate}
Applying L\^e- Greuel type formula \eqref{formulaGreuel} of 
Proposition \ref{LeGreuel} for the surfaces $X\cap H$ and $X\cap H'$, we have:
\begin{equation}
\label{eq:minimo1}
\mu (X\cap H)+ \mu(X\cap H\cap H')= m_2(X\cap H)
\end{equation}
\begin{equation}
\label{eq:minimo2}
\mu (X\cap H')+ \mu(X\cap H\cap H')= m_2(X\cap H')
\end{equation}
As $H$ and $H'$ are general to $X$ at  0, it follows by Proposition  \ref{Tessier} that $m_2(X\cap H)=m_2(X\cap H')= m_2(X)$.
Then, equations \eqref{eq:minimo1} and \eqref{eq:minimo2} imply that $\mu (X\cap H)=\mu (X\cap H')$, that is $\mu(X\cap H)$ is minimum.

Let $H'$ and $H''$ be  general hyperplanes to $X$ at the origin, defined by the kernels of $p'$ and $p''$, such that $\mu(X\cap H'\cap H'')$ is minimum. Then, the following L\^e- Greuel type formulas holds for $X\cap H$ and $X\cap H'$.
\begin{eqnarray}
\mu (X\cap H)+ \mu(X\cap H\cap H') & = & m_2(X\cap H)=m_2(X) \label{eq:minimo3} \\
\mu (X\cap H')+ \mu(X\cap H'\cap H'') & = & m_2(X\cap H')=m_2(X) \label{eq:minimo4}
\end{eqnarray}

As $\mu (X\cap H)=\mu (X\cap H')$ it follows that $\mu (X\cap H\cap H')=\mu (X\cap H'\cap H'')$.

Conversely, let $\pi=(p,p',p'')$ be the map $\C^N\to \C^3$, with $H=\ker p$, $H'=\ker p'$ and $H''=\ker p''$, such that $\ker \pi \cap C_{X,0}=\{0\}$ and $m_0(X)=\text{deg } \pi|_X$. Consider $\pi':\C^3\to \C^2$ such that $\pi' \circ \pi |_{X\cap H}=(p',p'')$, $\ker (\pi' \circ \pi) \cap C _{X\cap H}=\{0\}$ and $m_0(X)=m_0(X\cap H)= \text{deg } \pi' \circ \pi|_{X\cap H}$. Then
\begin{equation}
\label{ }
\mu(X\cap H \cap H')+\text{deg } \pi-1=m_1(X\cap H\cap H')
\end{equation}
\begin{equation}
\label{ }
\mu(X\cap H' \cap H'')+\text{deg } \pi-1=m_1(X\cap H'\cap H'').
\end{equation}
We have that $m_1(X)=m_1(X\cap H\cap H')$. On the other hand, it follows from  \cite[Lema 4.3]{Jawad} that
$$\mu(X\cap H \cap H')+\text{deg } \pi-1=(\Delta _{\pi' \circ \pi}\cdot (\pi' \circ \pi)(X\cap H\cap H'))_0,$$ where the notation $(a\cdot b)_0$ indicates the intersection multiplicity of $a$ and $b$. Then $m_1(X\cap H\cap H')=(\Delta _{\pi' \circ \pi}\cdot (\pi' \circ \pi)(X\cap H\cap H'))_0$  implies that $(\pi' \circ \pi)(X\cap H\cap H')$ is not limit of tangent hyperplanes to $\Delta _{\pi' \circ \pi}$.
Then by Proposition \ref{Legeral}, we have that $\pi(X\cap H) $ is not  limit of tangent hyperplanes to $\Delta_{\pi}$. It follows from Proposition \ref{Legeral} that $H$ is not  limit of tangent hyperplanes to $X$.
\end{prova}

\begin{exe}
 Let $X\subset \C^5$ be a 3-determinantal variety with isolated singularity defined by
$$\begin{array}{cccl}
      F:&\C^5& \to  & M_{2,3} \\
      &   (x,y,z,w,v)&\mapsto& \begin{pmatrix}
      x&y&z    \\
      w& v& x^2+y^2
\end{pmatrix}
\end{array}$$

Let $H$ and $H'$ be hyperplanes  given by the kernels of $p(x,y,z,w,v)=w-z$ and $p'(x,y,z,w,v)=x-v$. The  surfaces $X\cap H$ and  $X\cap H'$ are represented by the follo\-wing matrices  
$$\begin{pmatrix}
    x  &y & z    \\
     z & v & x^2+y^2
\end{pmatrix}  \text{  and  } \begin{pmatrix}
    x  &y & z    \\
     w& x & y^2
\end{pmatrix} \text{ respectively}.$$

It follows  (see  \cite{MIRIAMSP}) that $\mu(X\cap H)=4$ and $\mu(X\cap H')=2$, and that $\mu(X\cap H')=2$ is the minimum Milnor number.
\end{exe}

\section{Sections of EIDS}
In this section, we show that when $X$ is an EIDS of dimension $d>2$ it is possible to associate to $X$ a determinantal surface $Y$ whose Milnor number is an invariant of $X$.

\begin{prop}  \label{Le} \cite{Le}  Let $H$ and $H'$  be hyperplanes in $\C^N$ that are not limits of  tangent hyperplanes  to  $(X,0)$. Then there exist planes $Q$ and $Q'$ in $H$ and $H'$ respectively such that the reduced curves $|X\cap Q|$ and $|X \cap Q' |$ have the same Milnor number.
\end{prop}

 The following result  is a generalization of a result of L\^e D. T.
\begin{teo}  Let $X^d\subset \C^N$ be an EIDS and let $H$, $H'$ be strongly general hyperplanes  to $(X,0)$ at the origin. Then $H$ and $H'$ contain $P$ and $P'$ respectively  such that $\codim P= \codim P'=d-2$,  for which the determinantal surfaces $X\cap P$ and $X\cap P'$ satisfy  the following conditions:
\begin{itemize}
\item [a)] $X\cap P$ and $X\cap P'$ have isolated singularity.
\item [b)] $X\cap P$ and $X\cap P'$ admit smoothing.
\item [c)] $\mu(X\cap P)=\mu(X\cap P').$
\end{itemize}
\end{teo}

\begin{prova}
\begin{itemize}
 \item [a)] We can consider $P=H\cap \dots\cap H_{d-3}\subset H$ and $P'=H'\cap \dots \cap H'_{d-3}\subset H'$ such that, for all $i$, $H_i$ and $H'_i$ are strongly general for $X\cap H\cap H_1\cap \dots\cap H_{i-1}$ and $X\cap H' \cap H' _1\cap \dots\cap H' _{i-1}$ respectively.  By proposition \ref{EIDS}, it follows that $X\cap P\subset \C^{N-d+2}$ and $X\cap P'\subset \C^{N-d+2}$ are EIDS. As $\codim X< (m-t+2)(n-t+2)-2$ then $X\cap P$ and $X\cap P'$ are determinantal surfaces with isolated singularity.

\item [b)] Let $\mathcal{X}$ be an essential smoothing of  $X$, such that $H_i\times \C$ and $H'_i\times \C$  are strongly general hyperplanes for $\mathcal{X}\cap( H\times \C)\cap( H_1\times\C)\cap \dots\cap (H_{i-1}\times\C)$ and $\mathcal{X}\cap (H'\times \C) \cap (H' _1\times\C)\cap \dots(\cap H' _{i-1}\times \C)$ respectively. In fact, $\mathcal{X}\cap \widetilde{P}$ and $\mathcal{X}\cap \widetilde{P}'$, with $\widetilde{P}=P\times \C$ and $\widetilde{P}'=P'\times \C$ are the smoothings of $X\cap P$ and $X\cap P'$ respectively.

 \item [c)] The hyperplanes $H$ and $H'$ contain planes $Q=P\cap H'$ and $Q'= P' \cap H$ given in Proposition \ref{Le}. Then $X\cap Q$ and $X\cap Q'$ are determinantal curves  and by Proposition  \ref{Le}, we have that:
\begin{equation}\label{eq:curva}
\mu(X\cap Q)=\mu(X\cap Q').
\end{equation}
Using L\^e-Greuel type formula, Proposition \ref{LeGreuel}, we get
\begin{eqnarray}
\mu(X\cap P)+\mu(X\cap Q) & = & m_2(X\cap P, p') \\
\mu(X\cap P')+\mu(X\cap Q') & = & m_2(X\cap P', p).
\end{eqnarray}
We can apply Proposition \ref{Tessier} 
to the variety $\mathcal{X}\subset \C^N\times \C$, and to the hyperplanes $H_i\times \C$.
By definition, $m_2(X\cap P, p')=m_2(X\cap P', p)=m_0(P_2(\mathcal{X}))$, and $m_2(X\cap P', p)=m_0(P_2(\mathcal{X})).$

Then $ m_2(X\cap P, p')=m_2(X\cap P', p)=m_0(P_2(\mathcal{X}))$ and using  \eqref{eq:curva} we have the result.
\end{itemize}
\end{prova}
\begin{exe} Let $F:\C^N \to M_{2,3}$ be an analytic map, $N\geqslant 7$ defined by
 $$F(x,y)=\begin{pmatrix}
     x_1 &x_2 &x_3   \\
      x_4& x_5 & x_6+g(y)
\end{pmatrix}, $$
where $x=(x_1,x_2,\dots, x_6)$, $y=(y_1,\dots, y_{N-6})$ and $g:\C^{N-6}\to \C$ is an analytic function with $g(0)=0$. Then $X=F^{-1}(M_{m,n}^t)$ is a Cohen-Macaulay codimension 2 singularity in $\C^N$.
Let $P=\{(x,y)| \text{ }x_1=x_5, \text{ } x_2=x_6, \text{ } y=0\}$. Then $X\cap P$ is a determinantal surface in $\C^4$ defined by
$$\begin{pmatrix}
     x_1 &x_2& x_3   \\
   x_4   &  x_1& x_2
\end{pmatrix}$$
and  $\mu(X\cap P)=1$ (see \cite{MP}).

\end{exe}



\end{document}